\documentclass[10pt,twoside]{article}
\usepackage{graphicx}
\usepackage{amsmath}
\usepackage{amsfonts}
\usepackage{Latex-document}

\def\al{\alpha}

\def\om{\omega}
\def\la{\lambda}

\def\del{\overline\partial}
\def\ti{\times}
\def\M{{\cal M}}

\def\y{{\cal Y}}
\def\oy{\ov\y}
\def\Si{\Sigma}
\def\ov#1{\overline{#1}}

\def\Z{{ \mathbb Z}}
\def\R{{ \mathbb R}}
\def\P{{ \mathbb P}}
\def\Q{{ \mathbb Q}}
\def\cx{{ \mathbb C}}

\def\ra{\rightarrow}

\newtheorem{theorem}{Theorem}[section]
\newtheorem{prop}[theorem]{Proposition}

\newtheorem{conj}[theorem]{Conjecture}

\title{\bf  Symplectic Sums  and \vskip -2mm Gromov-Witten Invariants
\vskip 6mm}

\author{Eleny-Nicoleta Ionel\vspace*{-0.5cm}\thanks{
Department of Mathematics, University of Wisconsin-Madison,
Madison, WI 53706, USA. E-mail: ionel@math.wisc.edu}}
\date{\vspace{-8mm}}

\begin{document}

\maketitle

\thispagestyle{first} \setcounter{page}{427}

\begin{abstract}

\vskip 3mm

Gromov-Witten invariants of a symplectic manifold are a count of
holomorphic curves. We describe a formula
expressing the GW invariants of a symplectic sum $X\# Y$ in terms of
the relative GW invariants of $X$ and $Y$. This formula has several
applications to enumerative geometry. As one application, we obtain new
relations in the cohomology ring of the moduli space of complex
structures on a genus g Riemann surface with n marked
points.

\vskip 4.5mm


\end{abstract}

\vskip 12mm
\section{Gromov-Witten invariants} \label{section 1}\setzero

\vskip-5mm \hspace{5mm}

A symplectic structure on a closed smooth manifold $X^{2N}$ consists
of a closed, non-degenerate 2-form $\om$. Gromov's idea \cite{g} was
that one could obtain information about the symplectic structure on
$X$ by studying holomorphic curves. For that one needs to introduce an
almost complex structure, which is an endomorphism $J\in End(TX)$ with
$J^2=-Id$. Such a $J$ is compatible with $\omega$ if the bilinear form
$g(v,w)=\om(v,Jw)$ defines a Riemannian metric on $TX$. For a fixed
symplectic structure, the space of compatible almost complex
structures is a nonempty, contractible space.

One then considers the moduli space of $J$-holomorphic maps from
Riemann surfaces into $X$. Constraints are  imposed on the maps,
requiring the domain to have a certain form and the image to pass
through geometric representatives of fixed homology classes in
$X$. When the right number of constraints are chosen there will be finitely
many maps satisfying those constraints; the (oriented) count of these
maps will give the corresponding Gromov-Witten invariant. In general,
there are several technical difficulties one must overcome to get a
well-defined Gromov-Witten invariant. The foundations of this theory
began with \cite{g}, \cite{pw}, \cite{rt1} and have been developed since then
by the efforts of a large group of mathematicians (see, for example,
the references in \cite{ip3} and \cite{ms}). Here we present a brief
overview of the technical setup.

Consider $(X,\omega)$ a symplectic manifold. For each compatible
almost complex structure $J$ and perturbation $\nu$ one considers maps
$f:C\ra X$ from a genus $g$ Riemann surface $C$ with $n$ marked points
which satisfy the pseudo-holomorphic map equation $\del f=\nu$ and
represent a fixed homology class $A=[f]\in H_2(X)$. The set of such
maps (modulo reparametrizations), together with their limits, forms
the compact space of stable maps $\ov\M_{g,n}(X,A)$. For each stable
map $f:C\ra X$, the domain determines a point in the Deligne-Mumford
moduli space $\ov\M_{g,n}$ of genus $g$ Riemann surfaces with $n$
marked points (see also \S 3). The evaluation at each marked point
determines a point in $X$. All together, this gives a natural map
$$
\ov\M_{g,n}(X,A) \longrightarrow \ov\M_{g,n} \ti X^n.
$$
For  generic $(J,\nu)$ the image of this map carries a fundamental
homology class $[GW_{X,A,g,n}]$ which is defined to be the
Gromov-Witten invariant of $(X,\om)$. The dimension of this homology
class, given by an index computation, is
$$
\dim \ov\M_{g,n}(X,A)  =2c_1(TX)A+(\dim X-6)(1-g)+2n.
$$
A cobordism argument shows that the homology class $[GW_{X,A,g,n}]$ is
independent of generic $(J,\nu)$ and moreover depends only on the
isotopy class of the symplectic form $\omega$. Frequently, the
Gromov-Witten invariant is thought of as a collection of numbers
obtained by evaluating the homology class $[GW_{X,A,g,n}]$ on a basis
of the dual
cohomology group.  For complex algebraic manifolds these symplectic
invariants can also be defined by algebraic geometry, and in important
cases the invariants are the same as the counts of curves that are the
subject of classical enumerative algebraic geometry.

\medskip

The next important question is to find effective ways of computing the
GW invariants. One useful technique is the method of `splitting the
domain'. Anytime we have a relation in the cohomology of $\ov\M_{g,n}$
it pulls back to a relation (sometimes trivial) between the GW
invariants of a symplectic manifold $X$. As an example, suppose that
the constraints imposed on the domain of the holomorphic curves are
boundary classes in  $H^*(\ov\M_{g,n})$ (as defined in section 3
below). One then obtains recursive relations which relate such
GW invariant to invariants of lower degree or genus.  This method was
first used by Kontsevich and Ruan-Tian \cite{rt1} to determine
recursively the genus 0 invariants of the projective spaces $\P^n$. These
recursive relations follow from the observation that in the
Deligne-Mumford space $\ov\M_{0,4}\cong \P^1$ each boundary class
corresponds to a point, and are thus all homologous to each other.

In joint work with Thomas H. Parker, the author established a general
formula describing the behavior of GW invariants under the operation
of `splitting the target'  (\cite{ipa}, \cite{ip3}, \cite{ip4}).
Because we work in the context
of symplectic manifolds the natural splitting of the target is the one
associated with the symplectic cut operation and its inverse, the
symplectic sum. The next section describes the
symplectic sum operation and the main ingredients entering the sum
formula for GW invariants.

\bigskip

\section{Symplectic sums} \label{section 2} \setzero\vskip-5mm \hspace{5mm }

The operation of symplectic sum is defined by gluing along codimension
two submanifolds (see \cite{gf}, \cite{mw}).  Specifically, let $X$ be a
symplectic manifold with a codimension two symplectic submanifold $V$.
Given a similar pair $(Y,V)$ with a symplectic identification between
the two copies of $V$ and a complex anti-linear isomorphism between the
normal bundles $N_XV$ and $N_YV$ of $V$ in $X$ and in $Y$ we can form
the symplectic sum $X\#_{V}Y$.

Perhaps it is in more natural to describe the symplectic sum not as a single
manifold but as a family $Z\to D$ over the disk depending on a
parameter $\la\in D$.  For $\la\neq 0$ the  fibers $Z_\la$
are smooth and symplectically isotopic to $X\#_VY$ while the central
fiber $Z_0$ is the singular manifold $X\cup_V Y$. In a neighborhood of
$V$  the total space $Z$ is $N_XV\oplus N_YV$ and the fiber $Z_\la$
is defined by the equation $x y=\la$ where $x$ and $y$ are coordinates
in the normal bundles $N_XV$ and $N_YV\cong (N_XV)^*$.  The fibration $Z\to
D$ extends away from $V$ as the disjoint union of $X\times D$ and
$Y\ti D$.

Our overall strategy for proving the symplectic sum formula for GW
invariants \cite{ip4} is to relate the pseudo-holomorphic maps into
$Z_\la$ for $\la$ small to pseudo-holomorphic maps into $Z_0$. One
expects the stable maps into the sum to be pairs of stable maps into
the two sides which match in the middle. A sum formula thus requires a
count of stable maps in $X$  that keeps track of how the curves
intersect $V$.

So the first step is to construct Gromov-Witten invariants for a
symplectic manifold $(X,\omega)$ relative to a codimension two
symplectic submanifold $V$.  These invariants were introduced in
a separate paper with Thomas H. Parker \cite{ip3} and were designed
for use in symplectic sum formulas. Of course, before speaking of
stable maps one must extend the almost complex structure $J$ and the
perturbation $\nu$ to the symplectic sum.  To ensure that there is
such an extension we require that the pair $(J,\nu)$ be
$V$-compatible. The precise definition is given in section \S 6 of
\cite{ip3}, but in particular for such pairs $V$ is a $J$-holomorphic
submanifold --- something which is not true for generic $J$. The
relative invariant gives counts of stable maps for these special
$V$-compatible pairs.  Such counts are in general different from those
associated with the absolute GW invariants described in the first
section of this note.

Restricting to $V$-compatible pairs has repercussions. Any
pseudo-holomorphic map $f:C\to V$ into $V$ then automatically
satisfies the pseudo-holomorphic map equation into $X$. So for
$V$-compatible $(J,\nu)$, stable maps may have domain components whose
image lies entirely in $V$, so they are far from being transverse to
$V$.  Worse, the moduli spaces of such maps can
have dimension {\it larger} than the dimension of $\M_{g,n}(X,A)$.  We
circumvent these difficulties by restricting attention to the stable
maps which {\it have no components mapped entirely into $V$}.  Such
`$V$-regular' maps intersect $V$ in a finite set of points with
multiplicity.  After numbering these points, the space of $V$-regular
maps separates into components labeled by vectors $s=(s_1, \dots ,
s_{\ell})$, where $\ell$ is the number of intersection points and
$s_k$ is the multiplicity of the $k^{th}$ intersection point. Each
(irreducible) component $\M^{V}_{g,n,s}(X,A)$ of $V$-regular stable
maps is an orbifold; its dimension depends of $g,n,A$ and on the
vector of multiplicities $s$.

Next key step is to show  that the space of
$V$-regular maps carries a fundamental homology class.  For this we
construct an orbifold compactification $\ov \M^V_{g,n,s}(X,A)$, the space of
$V$-stable maps.  The relative invariants are then defined in exactly
the same way as the GW invariants. We consider  the natural map
\begin{equation}
\label{rel.GW} \ov\M^V_{g,n,s}(X,A)\ra  \ov \M_{g,n+\ell} \ti
X^n\ti V^\ell.
\end{equation}
The new feature is the last factor (the evaluation at the $\ell$
points of contact with $V$) which allows us to constrain how the
images of the maps intersect $V$.  Thus the relative invariants give
counts of $V$-stable maps with constraints on the complex structure of
the domain, the images of the marked points, and the geometry of the
intersection with $V$. There is one more complication: to be useful for a
symplectic sum formula, the relative invariant should record the
homology class of the curve in $X\setminus V$ rather than in $X$. This
requires keeping track of some additional homology data which is
intertwined with the intersection data, as explained in \cite{ip3}.

\medskip

We now return to the discussion of the symplectic sum formula. As
previously mentioned, the overall strategy is to relate the
pseudo-holomorphic maps into $Z_0$, which are simply maps into $X$ and $Y$
which match along $V$, with pseudo-holomorphic maps into $Z_\la$ for
$\la$ close to zero. For that we consider sequences of stable
maps into the family $Z_\la$ of symplectic sums as the `neck size'
$\la\to 0$.  These limit to maps into the singular manifold
$Z_0=X\cup_V Y$. A more careful look reveals several features of the
limit maps.

First of all, if the limit map $f_0:C_0\ra Z_0$ has no components in
$V$ then $f_0$ has matching intersection with $V$ on $X$ and $Y$
side. For such a limit map $f_0$ all its intersection points with $V$
are nodes of the domain $C_0$. Ordering this nodes we obtain a
sequence of multiplicities $s=(s_1,\dots, s_\ell)$ along $V$. But it
turns out that the squeezing process is not injective in general. For
a fixed $\la\ne 0$ there are $|s|=s_1\cdot \ldots \cdot s_\ell$ many
stable maps into $Z_\la$ close to $f_0$.

Second, connected curves in $Z_\la$ can limit to curves whose
restrictions to $X$ and $Y$ are not connected.  For
that reason the GW invariant, which counts stable curves from a
connected domain, is not the appropriate invariant for expressing a
sum formula.  Instead one should work with the `Gromov-Taubes'
invariant $GT$, which counts stable maps from domains that need not be
connected.
Thus we seek a formula of the general form
\begin{equation}
 GT_{X\#_V Y}\;=\; GT_X^{V}\, *\, GT_Y^{V}
\label{0.1}
\end{equation}
where $*$ is the operation that adds up the ways curves on the $X$ and
$Y$ sides match and are identified with curves in $Z_\la$.  That
necessarily involves keeping track of the multiplicities $s$ and the
homology classes. It also involves accounting for the limit maps
which have components in $V$; such maps are not counted by the
relative invariant and hence do not contribute to the left side of
(\ref{0.1}).

Finally, we need to consider limit maps which have components mapped
entirely in $V$. We deal with that possibility by squeezing the neck
not in one region, but several regions.  As a result, the formula
(\ref{0.1}) in general has an extra term called the $S$-matrix which
keeps track of how the genus, homology class, and intersection points
with $V$ change as the images of stable maps pass through the neck
region.  One sees these quantities changing abruptly as the map passes
through the neck --- the maps are ``scattered'' by the neck.  The
scattering occurs when some of the stable maps contributing to the GT
invariant of $Z_\la$ have components that lie entirely in $V$ in the
limit as $\la\to 0$.  Those maps are not $V$-regular, so are not
counted in the relative invariants of $X$ or $Y$.  But this
complication can be analyzed and related to the relative invariants of
the ruled manifold $\P(N_X V\oplus \cx)$.

Putting all these ingredients together, we can at last state the main
result of \cite{ip4}.
\begin{theorem}\label{GW}
Let $Z$ be the symplectic
sum of $(X,V)$ and  $(Y,V)$ and fix a decomposition of the constraints
$\al$ into $\al_X$ on the $X$ side and $\al_Y$ on the $Y$ side.
Then the GT invariant of $Z$ is given in terms of the relative
invariants of $(X,V)$ and $(Y,V)$  by
\begin{equation}
GT_Z(\al)\ =\   GT_X^{V}(\al_X) \, * \,
S_V\, * \, GT_Y^{V}(\al_Y)
\label{0.3}
\end{equation}
where $*$ is the convolution operation and $S_V$ is the $S$-matrix
defined in \cite{ip4}.
\end{theorem}
\medskip

Several applications of this formula are described in the next two
sections (see also \cite{ip4} for more applications). But the full
strength of the symplectic sum theorem has not yet been used.

A.-M. Li and Y. Ruan also have a sum formula \cite{lr}.  Eliashberg,
Givental, and Hofer are developing a general theory for invariants of
symplectic manifolds glued along contact boundaries \cite{egh}. Jun Li
has recently adapted our proof to the algebraic case \cite{li}.


\section{Relations in {\boldmath $H^*(\M_{g,n})$}} \label{section 3}
\setzero\vskip-5mm \hspace{5mm }

A smooth genus $g$ curve with $n$ marked points is stable if
$2g-2+n>0$. The set of such curves, modulo
diffeomorphisms, forms the moduli space $\M_{g,n}$. The stability
condition assures that the group of diffeomorphisms acts with finite
stabilizers, and so $\M_{g,n}$ has a natural orbifold structure. Its
Deligne-Mumford  compactification $\ov\M_{g,n}$ is a projective
variety. Elements of $\ov\M_{g,n}$ are called stable
curves; these are connected unions of smooth stable components $C_i$
joined at $d$ double points with a total of $n$ marked points and
Euler characteristic $\chi=2-2g+d$. The compactification $\ov\M_{g,n}$
is also an orbifold, and in fact Looijenga proved that it has a finite
degree cover which is a smooth manifold. In any event, the rational
cohomology of $\ov\M_{g,n}$ satisfies Poincar\'{e} duality. Throughout
this section we work only with rational coefficients.

There are several maps between moduli spaces of stable
curves. First, there is a projection $\pi_i:\ov\M_{g,n+1}\ra
\ov\M_{g,n}$ that forgets the marked point $x_i$ (and collapses the
components that become unstable). Second, we can consider the
attaching maps that build a boundary stratum in $\ov\M_{g,n}$. For
each topological type of a stable curve with $d$ nodes, with
components $C_i$ of genus $g_i$ and $n_i$ marked points the attaching
map $\xi$ at the $d$ nodes takes $\sqcup_i \ov\M_{g_i,n_i}$ onto a
boundary stratum of $\ov\M_{g,n}$.

We focus next on three kinds of natural classes in $H^*(\ov\M_{g,n})$
(or the Chow ring). For each $i$ between 1 and $n$ let $L_i\ra
\ov\M_{g,n}$ denote the relative cotangent bundle to the stable curve
at the marked point $x_i$.  The fiber of $L_i$ over a point
$C=(\Si,x_1,\dots,x_n)\in\ov\M_{g,n}$ is the cotangent space to $\Si$
at $x_i$, and its first Chern class $\psi_i$ is called a
{\em descendant} class. So there are $n$ descendant
classes $\psi_1, \dots, \psi_n$, one for each marked point. Next,
there are {\em tautological} (or Mumford-Morita-Miller) classes
$\kappa_0,\kappa_1,\dots $ obtained from powers of descendants by
the formula $\kappa_a=(\pi_{n+1})_*(\psi_{n+1}^{a+1})$ for each $a\ge
0$ (where $\pi_*$ denotes the push forward map in cohomology defined
using the Poincar\'{e} duality). Finally, the Poincar\'{e} dual of a
boundary stratum is called a {\em boundary} class. These three kinds of
natural classes are all algebraic and even dimensional; we define their
{\em degree} to be their complex dimension.

One natural --- and difficult --- problem is to describe the structure
of the cohomology rings of $\M_{g,n}$ and $\ov\M_{g,n}$.  This arises
from a different perspective as well since $H^*(\M_{g,n})$ is also the
cohomology of the mapping class group (for more details, see Tillman's
I.C.M. talk).  In genus zero Keel \cite{keel} determined the
cohomology ring of $\ov\M_{0,n}$ in terms of generators (which are
boundary classes) and relations. For higher genus far less is known
about the cohomology ring.

In this section  we will instead focus on finding relations in the
cohomology ring. For example, in genus 0 all relations come from the
``4-point relation'', essentially that in the
cohomology of
$\ov\M_{0,4}\cong \P^1$ the four $\psi_i$ classes as well as the three
boundary classes are all cohomologous (all being Poincar\'{e} dual to a
point). In genus 1 it is also known that $\psi_1$ is equal to 1/12 of
the boundary class in $\ov\M_{1,1}$. One might wonder whether in
higher genus all the $\psi$ classes come from the
boundary. That turns out not to be true in genus $g\ge 2$, but in
genus 2 Mumford \cite{mu} found a relation in $\ov\M_{2,1}$ expressing
$\psi^2_1$ as a combination of boundary classes. Several years
ago, Getzler
\cite{gz} found a similar relation for $\psi_1\psi_2$ in $\ov\M_{2,2}$
and he conjectured that this pattern would continue in higher
genus. In fact,
\begin{theorem}
\label{T.main}
When $g\ge 1$, any product of descendant or tautological classes of
degree at least $g$ (or at least $g-1$ when $n=0$) vanishes when
restricted to $H^*(\M_{g,n},\Q)$.
\end{theorem}

This result was proved by the author in \cite{i}. It extends an
earlier result of Looijenga \cite{l1}, who proved that a product
of descendant classes of degree at least $g+n-1$ vanishes in the
Chow ring $A^*({\cal C}^n_g)$ of the moduli space ${\cal C}^n_g$
of smooth genus $g$ curves with $n$ not necessarily distinct
points.

The idea of proof of Theorem \ref{T.main} is simple. We start
with the moduli space $\y_{d,g,n}$ of degree $d$ holomorphic maps from
smooth genus $g$ curves with $n$ marked points to $S^2$ which have a
fixed ramification pattern over $r$ marked points in the target.  We
then consider its relative stable map compactification $\ov\y_{d,g,n}$
(closely related to the space of admissible covers \cite{hm}). The
space $\ov\y_{d,g,n}$ has an orbispace structure and it comes with two
natural maps $st$ and $q$ that record respectively the domain and the
target of the cover.
\begin{equation}
\begin{array}{lcr}
&\oy_{d,g,n}&
\\
\hskip.1in ^{st}\swarrow&&\searrow^q\hskip.2in
\\
\ov \M_{g,n}&&\ov\M_{0,r}
\end{array}
\end{equation}
A simple way to get relations in the cohomology of $\ov\M_{g,n}$ is to
  pull back by $q$ known relations in the cohomology of $\ov\M_{0,r}$,
and then push them forward by $st$.

To begin with, note that the diagram above provides several other
natural classes in $\ov\M_{g,n}$: for each choice of ramification
pattern, $st_*\oy_{d,g,n}$ defines a cycle in $\ov\M_{g,n}$. The most
useful ones turn out to be the ``2-point ramification cycles'', for
which all but at most two of the branch points are simple. Pushing
forward such cycles by the attaching map of a boundary stratum gives a
generalized 2-point cycle.

To prove Theorem \ref{T.main}, we choose a degree $d$ of the cover
and a 2-point ramification cycle $\oy_{d,g,n}$ in such a way that  the
stabilization map $st:\oy_{d,g,n}\ra \ov\M_{g,n}$ has finite, nonzero
degree. The key step is the following proposition.
\begin{prop}\label{prop.PD} The Poincar\'{e} dual of any degree $m$
product of descendant and tautological classes can be written as a linear
combination of generalized 2-point ramification cycles of codimension
$m$.
\end{prop}

But the codimension of a 2-point ramification cycle is at most
$g$. A simple degeneration argument proves that the cycles of
codimension exactly $g$ vanish on $\M_{g,n}$, thus implying
Theorem \ref{T.main}.

There are three main ingredients in the proof of Proposition
\ref{prop.PD}. First, the relative cotangent bundle to the domain is
related to the pullback of the relative cotangent bundle to the
target, so we can express the descendant classes in the domain via
descendant classes in the target. Second, the target has genus zero
and (nontrivial) products of descendants in $\ov\M_{0,r}$ are Poincar\'{e}
dual to boundary cycles $D$. This means that we can relate a product of
descendants on the domain to cycles of type $st_*q^* D$. Finally,
a degeneration formula, which is essentially a consequence of
the symplectic sum Theorem \ref{GW}, expresses cycles of type
$st_*q^* D$ in terms of 2-point ramification cycles.

\smallskip

The degree $g$ in Theorem \ref{T.main} is the lowest degree in which
some {\em monomial} in
descendants would vanish on $\M_{g,n}$ (see the discussion in
\cite{i2}). However, there are lower degree polynomial relations in
descendent and tautological classes. For example, if we restrict our
attention to the moduli space $\M_g$ of smooth genus $g$ curves then
the subring generated by the tautological classes is called the {\em
tautological ring} $R^*_g$. Looijenga's result \cite{l1}
implies that $R^*_g=0$ for $*\ge g-1$ and  Faber \cite{fa} made the
following
\begin{conj} The classes $\kappa_1,\dots,\kappa_{[g/3]}$ generate the
tautological ring $R^*_g$.
\end{conj}

We refer the reader to \cite{fa} for the full conjecture.
\smallskip

It turns out that techniques similar to those of Theorem  \ref{T.main}
produce several other sets of relations between
tautological classes. One such set of relations implies
that, for each $a>[g/3]$, the class $\kappa_a$ can be written as
polynomial in lower degree tautological classes, as required by
Faber's conjecture. A detailed proof will appear in \cite{i}.

\section{Further applications} \label{section 4} \setzero\vskip-5mm
\hspace{5mm }

There are other applications of the sum formula (\ref{0.3}). One such
application considered in \cite{ip4} begins with the following simple
observation. Given any symplectic manifold $X$ with a codimension 2
symplectic submanifold $V$, we can write $X$ as a (trivial) symplectic
sum $X\#_V P_V$ where $P_V$ is the ruled manifold $\P(N_XV\oplus \cx)$
and $V$ is identified with its infinity section. We can then obtain
recursive formulas for the GW invariants of $X$ by moving constraints
from one side to the other and applying the symplectic sum formula.

In \cite{ip3} we used this method to obtain both (a) the
Caporaso-Harris formula for the number of nodal curves in $\P^2$
\cite{ch}, and (b) the ``quasimodular form'' expression for the
rational enumerative invariants of the rational elliptic surface
\cite{bl}.  In hindsight, our proof of (a) is essentially the same
as that in \cite{ch}; using the symplectic sum formula makes the proof
considerably shorter and more transparent, but the key ideas are the
same.  Our proof of (b), however, is completely different from that of
Bryan and Leung in \cite{bl}.

\smallskip
We end  with another interesting application of the
Symplectic Sum Theorem \ref{GW}. For each  symplectomorphism $f$ of a
symplectic manifold $X$, one can form the {\em symplectic mapping cylinder}
\begin{equation}
X_f = {X\times \R\times S^1}/\Z
\label{eq00.1}
\end{equation}
where the $\Z$ action is generated by $(x,s,\theta)\mapsto
(f(x),s+1,\theta)$. In a joint paper \cite{ip2} with T. H. Parker
we regarded $X_f$ as a symplectic sum and computed the Gromov
invariants of the manifolds $X_f$ and of fiber sums of the $X_f$ with
other symplectic manifolds.  The result is a large set of interesting
non-K\"{a}hler symplectic manifolds with computational ways of
distinguishing them.  In dimension four this gives a symplectic
construction of the `exotic' elliptic surfaces of Fintushel and Stern
\cite{fs}.  In higher dimensions it gives many examples of manifolds
which are diffeomorphic but not `equivalent' as symplectic manifolds.

More precisely, fix a symplectomorphism $f$ of a closed symplectic
manifold $X$, and let $f_{*k}$ denote the induced map on
$H_k(X;\Q)$. Note that $X_f$ fibers over the torus $T^2$ with fiber
$X$.  If $\mbox{det }(I-f_{*1})=\pm 1$ then there is a well-defined
section class $T$.  Our main result of \cite{ip2} computes the genus
one Gromov
invariants of the multiples of this section class. These are the
particular GW invariants that, in dimension four,  C.H. Taubes
related to the Seiberg-Witten invariants(see \cite{t} and \cite{ip1}).

\begin{theorem}\label{IntroThm1}
 If $\mbox{det }(I-f_{*1})=\pm 1$, the partial Gromov series of $X_f$
for the section class $T$ is given by the Lefschetz zeta function of
$f$ in the variable $t=t_T$:
$$\displaystyle
Gr^T(X_f)\ =\ \zeta_f(t)\ =\ \frac{\prod_{k\  odd}\ \det(I-tf_{*k})}
{\prod_{k\  even}\
\det(I-tf_{*k})}.
$$
\end{theorem}

When  $X_f$ is a four-manifold,  a wealth of  examples arise from
knots.  Associated to each fibered knot $K$ in $S^3$ is a Riemann
surface $\Si$ and a monodromy diffeomorphism $f_K$ of $\Si$.
Taking $f=f_K$ gives symplectic 4-manifolds $X_K$ of the homology
type of  $S^2\ti T^2$ with
\[
Gr(X_K)= \frac{A_K(t_T)}{(1-t_T)^2}
\]
where $A_K(t)=\det(I-tf_{*1})$ is the Alexander polynomial of $K$ and $T$
is the section class.

 We can elaborate on this construction by fiber summing $X_f$ with
other 4-manifolds. For example, let $E(n)$ be the simply-connected
minimal elliptic surface with fiber $F$ and holomorphic Euler
characteristic $n$. Then $E(1)$ is the rational elliptic surface and
$K3=E(2)$.  Forming the fiber sum of $X_K$ with $E(n)$ along the tori
$T=F$, we obtain a symplectic manifold
$$
E(n,K)\ =\ E(n)\#_{F=T} X_K.
$$
homeomorphic to $E(n)$. In fact, for fibered knots $K,\; K'$ of the same
genus there is a homeomorphism between $E(n,K)$ and $E(n,K')$
preserving the periods of $\omega$ and the canonical class $\kappa$. For $n>1$
we can compute the full (not just partial) Gromov series.
\begin{prop}
\label{introprop2}
For $n\geq 2$, the Gromov and Seiberg-Witten series of $E(n,K)$ are
\begin{equation}\label{introThm2.1}
Gr(E(n,K))\ =\ SW(E(n,K))\ =\  A_K(t_F)\ (1-t_F)^{n-2}.
\end{equation}
\end{prop}

Thus fibered knots with distinct Alexander polynomials give rise
to symplectic manifolds $E(n,K)$ which are homeomorphic but not
diffeomorphic. In particular, there are infinitely many distinct
symplectic 4-manifolds homeomorphic to $E(n)$. Fintushel and Stern
\cite{fs} have independently shown how (\ref{introThm2.1}) follows
from knot theory and results in  Seiberg-Witten theory.

\label{lastpage}


\begin{thebibliography}{aa}

\bibitem{bl} J. Bryan  and N.-C. Leung, {\em The enumerative geometry
of $K3$ surfaces and modular forms},  J. Amer. Math. Soc. {\bf 13}
(2000), 371--410.

\bibitem{ch} L. Caporaso and  J. Harris, {\it Counting plane curves in
any genus}, Invent. Math. {\bf 131} (1998), 345--392.

\bibitem{egh} Y. Eliashberg, A. Givental and H. Hofer, {\em Introduction
to Symplectic Field Theory}, GAFA 2000 (Tel Aviv, 1999), Geom.
Funct. Anal. 2000, Special Volume, Part II, 560--673.

\bibitem{fa} C. Faber, {\em A conjectural description of the
tautological ring of the moduli space of curves}, Moduli of curves
and abelian varieties, 109--129, Aspects Math., E33, Vieweg,
Braunschweig, 1999.

\bibitem{fs} R. Fintushel and R. Stern, {\it Knots, Links and
4-Manifolds}, Invent. Math. {\bf 134} (1998), 363--400.

\bibitem{gz} E. Getzler, {\em Topological recursion relations in genus
2}, Integrable systems and algebraic geometry (Kobe/Kyoto, 1997),
73--106, World Sci. Publishing, River Edge, NJ, 1998.

\bibitem{gf} R. Gompf, {\it A new construction of symplectic
manifolds}, Annals of Math., {\bf 142} (1995), 527--595.

\bibitem{g} M. Gromov, {\it Pseudo holomorphic curves in symplectic
manifolds}, Invent. Math. {\bf 82} (1985), 307--347.


\bibitem{hm} J. Harris, I. Morrison, {\em Moduli of curves}, Graduate
Texts in Math, vol 187, Springer-Verlag, 1998.

\bibitem{i2} E. Ionel, {\it Topological recursive relations in
$H^{2g}({\cal M}_{g,n})$}, to appear in Invent. Math.

\bibitem{i} E. Ionel, {\it On relations in the tautological ring of
$\M_g$}, in preparation.

\bibitem{ip1} E. Ionel and T. H. Parker, {\it  The  Gromov invariants of
Ruan-Tian and Taubes},  Math. Res. Lett. {\bf 4} (1997), 521--532.

\bibitem{ip2} E. Ionel and T. H. Parker, {\it Gromov Invariants and
Symplectic Maps}, Math. Annalen, {\bf 314}, 127--158 (1999).

\bibitem{ipa} E. Ionel and T. H. Parker, {\it Gromov-Witten Invariants
of Symplectic Sums}, announcement, Math. Res. Lett., {\bf
5}(1998), 563--576.

\bibitem{ip3} E. Ionel and T. H. Parker, {\it Relative Gromov-Witten
Invariants}, to appear in Annals of Math.

\bibitem{ip4} E. Ionel and T. H.  Parker, {\it The Symplectic Sum Formula
for Gromov-Witten Invariants}, preprint, math.SG/0010217.

\bibitem{keel} S. Keel, {\em Intersection theory of moduli space of
stable $n$-pointed curves of genus zero}, Trans. Amer. Math. Soc.
{\bf 330}(1992),  545--574.

\bibitem{lr} A.-M. Li, Y. Ruan, {\it Symplectic surgery and
Gromov-Witten invariants of Calabi-Yau 3-folds}, Invent. Math.
{\bf 145} (2001), 151--218.

\bibitem{li} Jun Li, {\it A Degeneration formula of GW-invariants},
preprint, math.AG/0110113.

\bibitem{l1} E. Looijenga, {\em On the tautological ring of $\M_g$},
Invent. Math. {\bf 121}(1995), 411--419.

\bibitem{mw} J. McCarthy and J.Wolfson, {\it Symplectic Normal Connect
Sum}, Topology, {\bf 33} (1994) 729--764.

\bibitem{ms} D. McDuff and D. Salamon, {$J$-holomorphic curves and
quantum cohomology}, A.M.S., Providence, R.I., 1994.


\bibitem{mu} D. Mumford, {\em Towards an enumerative geometry of the
 moduli space of curves} in {\em Arithmetic and geometry II}
(ed. M. Artin and J. Tate), Progress in Math, vol 36, Birkh\"auser,
 Basel, 1983.

\bibitem{pw} T. H. Parker and J. Wolfson, {\it Pseudo-holomorphic maps
and bubble trees}, Jour. Geometric Analysis, {\bf 3} (1993)
63--98.

\bibitem{rt1} Y. Ruan and G. Tian, {\it A mathematical theory of
quantum cohomology}, J. Differential Geom. {\bf 42} (1995),
259--367.

\bibitem{rt2} Y. Ruan and G. Tian, {\it Higher genus symplectic
invariants and sigma models coupled with gravity}, Invent. Math.
{\bf 130} (1997), 455--516.


\bibitem{t} C. H. Taubes, {\it Counting pseudo-holomorphic curves in
dimension four}, J. Diff. Geom. {\bf 44} (1996), 818--893.

\end{thebibliography}
\end{document}